\documentclass[12pt]{article}

\usepackage{amssymb}
\usepackage{theorem}

\newtheorem{theorem}{Theorem}
\newtheorem{proposition}[theorem]{Proposition}
\newtheorem{corollary}[theorem]{Corollary}
\newtheorem{lemma}[theorem]{Lemma}
{\theorembodyfont{\rmfamily}
\newtheorem{definition}[theorem]{Definition}
\newtheorem{example}[theorem]{Example}

}

\def\cqfd{{$\diamond$}}
\def\X{{\bf X}}
\def\Aut{{\rm Aut }}
\def\R{{\mathbb{R}}}
\def\C{{\mathbb{C}}}
\def\CP{\mathbb{C}{\rm P}}
\def\tr{\, {\rm Tr}\,}

\include{epsf}

\date{\today}

\author{Dmitri Panov\thanks{panov@ihes.fr, \'Ecole Polytechnique,
91128 Palaiseau, France}, \,  
Dimitri Zvonkine\thanks{dimitri.zvonkine@math.u-psud.fr,
b\^at 425, Universit\'e Paris-Sud, 91400 Orsay, France}}

\title{Counting Meromorphic Functions with Critical
Points of Large Multiplicities}

\begin{document}

\maketitle

\begin{abstract}

We study the number of meromorphic functions on
a Riemann surface with given critical values
and prescribed multiplicities of critical points
and values.

When the Riemann surface is $\CP^1$ and the function
is a polynomial, we give an elementary way of finding
this number.

In the general case, we show that, as the multiplicities
of critical points tend to infinity, the asymptotic for the number
of meromorphic functions is given by the volume of some
space of graphs glued from circles. We express this
volume as a matrix integral.

\end{abstract}

\section{From functions to constellations}

In this section we sketch the more or less standard
construction that allows one 
to reduce the enumeration problem of meromorphic
functions to a combinatorial problem of enumerating
constellations. In the sequel we will be concerned with
the latter problem. The construction is based on the Riemann
existence theorem, exposed, for example, in~\cite{Volklein}.

Consider a compact connected Riemann surface 
$\Sigma$ and a non-constant meromorphic function $f$ on it.
The function $f$ can be considered as a ramified covering
of the Riemann sphere $\CP^1$ by the surface $\Sigma$.
Denote by $w_1, \dots, w_k \in \C$ the {\em finite} ramification
points of the covering. (The point $\infty \in \CP^1$ 
might also be a ramification point.) Let $w_0 \in \C$ be a finite 
non-critical value of $f$. Choose $k$ non-intersecting paths
on the plane, connecting $w_0$ to the points $w_1, \dots, w_k$.
The cyclic order of the paths at the point $w_0$ should
be equal to the order $w_1,\dots,w_k$.

\begin{figure}[htb] 
\begin{center}
\
\epsfbox{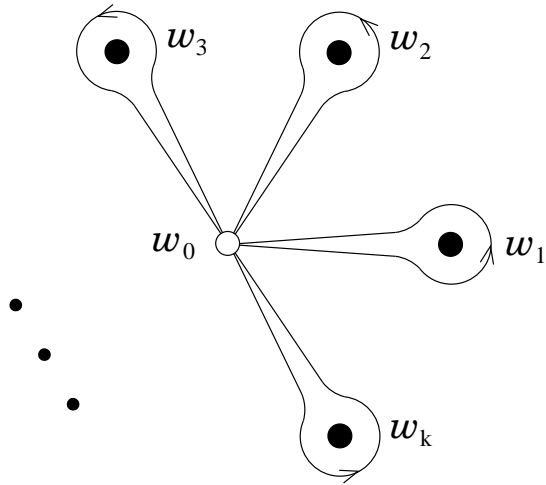}
\caption{A star of loops. \label{loopstar}}
\end{center}
\end{figure}

\begin{definition}
A {\em star of loops} on the complex plane is a set of $k$
non-intersecting closed loops $l_1, \dots l_k$ of the following form.
Each loop $l_i$ starts at $w_0$,
follows the path from $w_0$ to $w_i$, goes around $w_i$
in the counterclockwise direction, and then goes back to
$w_0$ following the path. A star of loops is shown
in Figure~\ref{loopstar}.
\end{definition}

In the sequel we suppose that the points $w_0,w_1, \dots, w_k$
and the star of loops are fixed once and for all.

The preimage under $f$ of the loop $l_i$ surrounding the critical
value $w_i$ is composed of several curvilinear polygons
on the Riemann surface $\Sigma$. The vertices of the
polygon are the preimages of $w_0$. Its sides 
are the preimages of the loop $l_i$. Each polygon surrounds a 
unique point $z$ such that $f(z) = w_i$. If the polygon has
more than one side, the corresponding point is critical;
if it has only one side, it is not. Now we take the preimages
of all the loops $l_i$ and erase all the 1-gons. We say
that the polygons that compose the preimage of $l_i$
are of {\em color} $i$.

Thus, on $\Sigma$, we obtain a graph consisting of 
polygons colored in $k$ colors, 
each polygon having at least 2 sides and
surrounding exactly one critical point of $f$. The vertices
of the polygons are the preimages of $w_0$. The polygons
meeting at a vertex are all of different colors, and the cyclic
order of the polygons meeting at a vertex is given by
the increasing order of colors.

\begin{definition}
A {\em $k$-constellation} is a connected colored
oriented graph consisting of
oriented polygons (with at least 2 sides each) colored in $k$ colors
(from $1$ to $k$) and such that no two polygons of the same
color share a common vertex.
\end{definition}

A constellation can be uniquely (up to a homeomorphism)
embedded into a (compact connected oriented) surface in such
a way that (i) if we cut the surface along the constellation
we obtain a union of domains homeomorphic to open discs,
(ii) every polygon of the constellation surrounds one such 
disc and is oriented counterclockwise with respect to this
disc, (iii) the cyclic order of polygons meeting at every
vertex coincides with the increasing order of colors.

\begin{definition}
The embedding satisfying the conditions (i), (ii), and (iii)
is called the {\em natural embedding} of a constellation.
\end{definition}

It follows from Riemann's existence theorem that (once the
star of loops is fixed) a 
constellation uniquely determines the couple $(\Sigma,f)$
up to a biholomorphic equivalence. The part of $\Sigma$ exterior to
the polygons consists of several connected components corresponding
to the poles of $f$. They will be called {\em faces}.

If $\Sigma$ is a sphere and $f$ a polynomial, the 
corresponding constellation is of a particular form, because
$\infty$ is always a critical value
of $f$ and its only preimage is $\infty$.
One can easily see that the constellation is, in that case,
a planar ``tree'' consisting of colored polygons.
Such a constellation is called a cactus.

\begin{definition}
A {\em $k$-cactus} is a constellation that has no cycles
other than the polygons it is glued from.
\end{definition}

The surface of the natural embedding of a cactus is
a sphere. Once a star of loops is fixed, a cactus
determines the polynomial $f$
uniquely up to a change of variable 
$f(z) \mapsto f(az+b)$.

\bigskip

Thus the problem of enumerating meromorphic functions
can be reduced to a combinatorial problem of enumerating
cacti and constellations, and we will concentrate on
the latter problem from now on.

\section{Spaces of gluings of circles}

The idea of this section is to replace by circles
the polygons
that compose cacti and constellations. The circles 
will be endowed with
elements of length. This corresponds to letting the 
numbers of sides of the polygons tend to infinity. The
number of cacti or constellations will be replaced
by the volume of the space of possible gluings of circles.
It turns out that these volumes are easier to compute then
the numbers of cacti and constellations.

\begin{definition} A {\em circle constellation}
is an oriented colored graph with nonnegative lengths
assigned to all edges, obtained by gluing together (at a 
finite number of points)
several circles of given lengths colored in $k$ colors
from $1$ to $k$.
Two circles of the same color cannot
share a common point.
\end{definition}

\begin{definition} A {\em circle cactus} is
a constellation of circles that has no cycles other 
than the circles themselves.
\end{definition}

As above, a circle constellation (or a circle cactus) has
a unique natural embedding into a compact surface (a
sphere), satisfying the
same three conditions as for ordinary constellations:
(i) if we cut the surface along the circles we obtain
a union of domains homeomorphic to open discs, (ii)
each circle of the circle constellation surrounds one
such disc and is oriented counterclockwise with respect
to it, (iii) the cyclic order of circles at every vertex
coincides with the increasing order of colors.

The connected components of the part of the surface
exterior to the circles are called {\em faces}.

Suppose a set of $m$ oriented colored 
circles with fixed lengths is given.
Let us fix a genus $g$ and a number of faces $p$. We can then
consider the space of all possible gluings of the circles that
give a circle constellation with $p$ faces on a surface of genus $g$.

\begin{proposition}
The set of circle constellations of genus $g$ with $p$ faces
glued of a given set of $m$ circles has a natural structure of
a singular smooth noncompact manifold of dimension $d=4g-4+m+2p$
with a volume measure. 
\end{proposition}

\paragraph{Proof.} First we will describe the space
of circle constellations with numbered vertices.
Such constellations do not have nontrivial 
automorphisms, which simplifies the description.

Let us describe the smooth part of the space of 
circle constellations. A circle constellation belongs to
the smooth part if all its circles intersect only two by two.
In other words, three circles never share a common vertex.
Let $\lambda_1, \dots, \lambda_q$ be the lengths of
the arcs into which one of the circles is divided by
the vertices. The sum of these lengths is equal to the length
of the circle (which is fixed), but, for example, the first
$q-1$ lengths can take arbitrary values in an open
domain in $\R^{q-1}$. Thus, taking for each circle 
the lengths of all its arcs except one, we obtain a
set of local coordinates of the space of circle constellations.
It is easy to check that the total number of these
arcs is indeed equal to $d$.

Now consider a smooth family of circle constellations depending on
a parameter $t$, such that, as $t$ tends to $0$, the lengths
of some of the arcs tend to $0$. Then the gluing that we 
obtain in the limit is still a gluing of the same circles
embedded into a surface of genus $g$ with $n$ faces, but
the cyclic order of the circles at the vertices is not
necessarily the same as the increasing order of the colors.
If the cyclic order of the circles at each vertex happens to
coincide with the increasing order of colors, then what we
have obtained is a true circle constellation. It is the
limit of the above family and lies in the nonsmooth
part of the space of gluings. If, on the contrary,
the cyclic order of circles at at least one vertex
differs from the increasing order of colors, then the
above family has no limit. This accounts for the
noncompactness of the space of gluings.

Now we can introduce the volume measure.
Consider the lengths of arcs $\lambda_i$, $1 \leq i \leq d$
that form a set of local coordinates on the smooth
part of the space of gluings. Consider the differential form 
$d\lambda_1 \wedge \dots \wedge d\lambda_d$. 
If we change the numbering of the
arcs, this differential form will either change
its sign or remain unchanged. Thus its absolute value
is a well-defined volume measure on the space of gluings.
We let the measure of the nonsmooth part be equal to~0.

If we now consider the space of circle constellations whose
vertices are not numbered, it is obtained from the above space
of circle constellations with numbered vertices by factoring
by a finite group action (the renumbering of the
vertices). Therefore it is also a smooth manifold
with singularities. The volume measure being indepenent
of the numbering of the vertices, it descends to the factor
space.
\cqfd

In the sequel we will also consider spaces
of circle constellations with a fixed number of marked points on the 
circles. An obvious modification of the above construction
shows that these spaces are also smooth manifolds with singularities
and that they are endowed with a natural volume form. To see
that, it suffices to cut every circle into arcs not only
at the points where it is glued to other circles, but also
at the marked points.

\bigskip

Now we are going to compare the number of constellations
composed of a given set of polygons and the volume of
the space of circle constellations composed of a given
set of circles.

Suppose we are given a set of polygons colored in $k$
colors and a set of circles colored in $k$ colors.
For each color, there is the same number of polygons
as of circles. Denote by $n_{ij}$ be the number of sides of the
$j$th polygon of color $i$ and by $l_{ij}$ the length
of the $j$th circle of color $i$. We assume that polygons of the
same color and number of sides are labeled, so that all the
polygons are distinguishable. Similarly, the circles of the
same color and lengths are labeled, so that all the
circles are distinguishable. 

Consider all constellations and circle constellations of genus
$g$ with $p$ faces. Denote by $P(n_{ij})$ 
(respectively $Q(l_{ij})$) the 
number of constellations (respectively the
volume of the space of circle constellations) like that
composed of the above set of polygons (circles). Finally,
let $m$ be the total number of polygons (circles) and
$d = 4g-4+m+2p$.

\begin{theorem} \label{asymptotic}
$P$ is a polynomial in the variables $n_{ij}$ of degree $d$.
$Q$ is a homogeneous polynomial in the variables $l_{ij}$ of degree $d$.
Moreover, $Q$ is obtained by taking the top degree homogeneous
part of $P$ and substituting $l_{ij}$ for $n_{ij}$.
In particular, if we substitutes $n_{ij}$ for
$L_{ij}$ in $Q(l_{ij})$, we obtain an asymptotic of the numbers
$P(n_{ij})$ for large $n_{ij}$.
\end{theorem}

\paragraph{Proof.} Suppose a constellation is given. We call the {\em
topological type} of the constellation the following information.
First, for each polygon we retain the number of vertices at which it
is glued to some other polygons. Second, enumerating these vertices
in the cyclic order in which they appear on the polygon, we list,
for each vertex, the polygons that are glued to it. The information
that we do not retain is the number of edges that lie between the listed
vertices.

It is easy to see that the number of constellations of a given
topological type is a polynomial in the $n_{ij}$. Indeed, suppose
that the $j$th polygon of color $i$ is glued to other
polygons at $k_{ij}$ vertices. Then there are
$$
{n_{ij}-1 \choose k_{ij}-1}
$$
ways to choose these vertices on the polygon. The total number
of constellations of a given topological type is equal to the
product of the above numbers over all the polygons (because
the choices of the vertices are independant for different
polygons). The top
degree terms are obtained from the topological types in which
no 3 polygons share a common vertex.

A topological type of a circle constellation is defined in exactly
the same way as for constellations: for each circle we list the points
at which it is glued to other circles and the circles that are glued at
each of these points; 
the lengths of arcs between these points are not retained.

If in a topological type at least 3 circles share a common point at least
once, than such a type determines a sub-manifold of the space of gluings of
positive codimension. Therefore its contribution to the total volume is 0.
If, on the other hand, no 3 circles share a common point, the topological
type determines a component of the space of gluings with nonzero volume.
As above, denote by $k_{ij}$ the number of points on the $j$th circle
of color $i$ at which it is glued to other circles. The volume of the
space of choices of these points equals
$$
\frac1{(k_{ij}-1)!} \, l_{ij}^{k_{ij}-1}.
$$
The volume of the whole component determined by the topological type
is the product of these numbers over all circles.

Since taking the top degree term of
$$
{n_{ij}-1 \choose k_{ij}-1}
$$
and replacing $n_{ij}$ by $l_{ij}$ gives us precisely
$$
\frac1{(k_{ij}-1)!} \, l_{ij}^{k_{ij}-1},
$$
the theorem follows. \cqfd

\begin{example}
Suppose we are given $3$ polygons (circles) of different
colors and with numbers of sides $n_1$, $n_2$, $n_3$
(respectively with lengths $l_1$, $l_2$, $l_3$). Then
we have $P = n_1+n_2+n_3-2$, $Q= l_1+l_2+l_3$
(see Examples~\ref{ex:3circles} and~\ref{ex:3polyg}).
\end{example}

\begin{theorem} \label{1circ}
The volume of the space of circle cacti glued from $k$ circles
of pairwise distinct colors and with lengths $l_1,\dots,l_k$
is equal to
$$
(l_1 + \dots + l_k)^{k-2}.
$$
\end{theorem}

\begin{example} \label{ex:3circles}
Let $k=3$. Then the three circles are necessarily glued
together in a chain. Suppose the circle number 1 is in the middle
(Fig.~\ref{3circles}). 
 
\begin{figure}[htb]
\begin{center}
\
\epsfbox{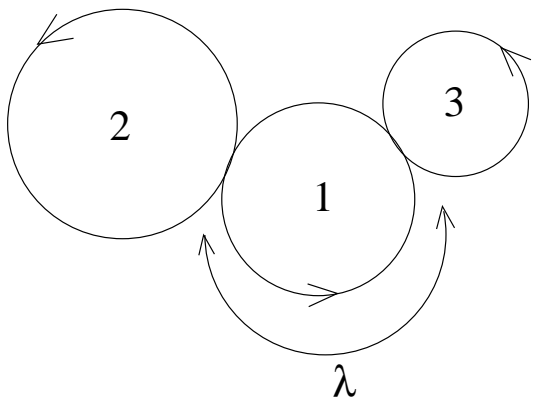}
\caption{A circle cactus with 3 circles. \label{3circles}}
\end{center}
\end{figure}

Then the length $\lambda$ of the arc between
the points at which it is glued to the circles 2 and 3 can
vary from $0$ to $l_1$. Thus the component of the space of 
gluings that corresponds to the first circle being in the
middle has volume $l_1$. Similarly, there are two more
components with volumes $l_2$ and $l_3$. The total volume
of the space of circle cacti is equal to $l_1 + l_2 + l_3$.
\end{example}

\paragraph{Proof of Theorem~\ref{1circ}.} 
We are going to describe a volume-preserving tranformation 
that allows one to reduce by $1$ the
number of circles (and of colors) by merging two circles into
one and adding their lengths. 
On the new circle cactus thus obtained one point will be marked,
so that the initial circle cactus can be recovered from the new one.
After applying this merging procedure $k-1$ times, we will
obtain a unique circle of length $l_1 + \dots + l_k$ 
with $k-1$ marked points. The space of possible markings
will thus have the same volume as the space of circle cacti,
and this volume is equal to
$$
(l_1 + \dots + l_k)^{k-2}.
$$

We start by making a list of the lengths $l_1, \dots, l_k$.
They will be needed to reverse the merging procedure.

\begin{figure}[htb]
\begin{center}
\
\epsfbox{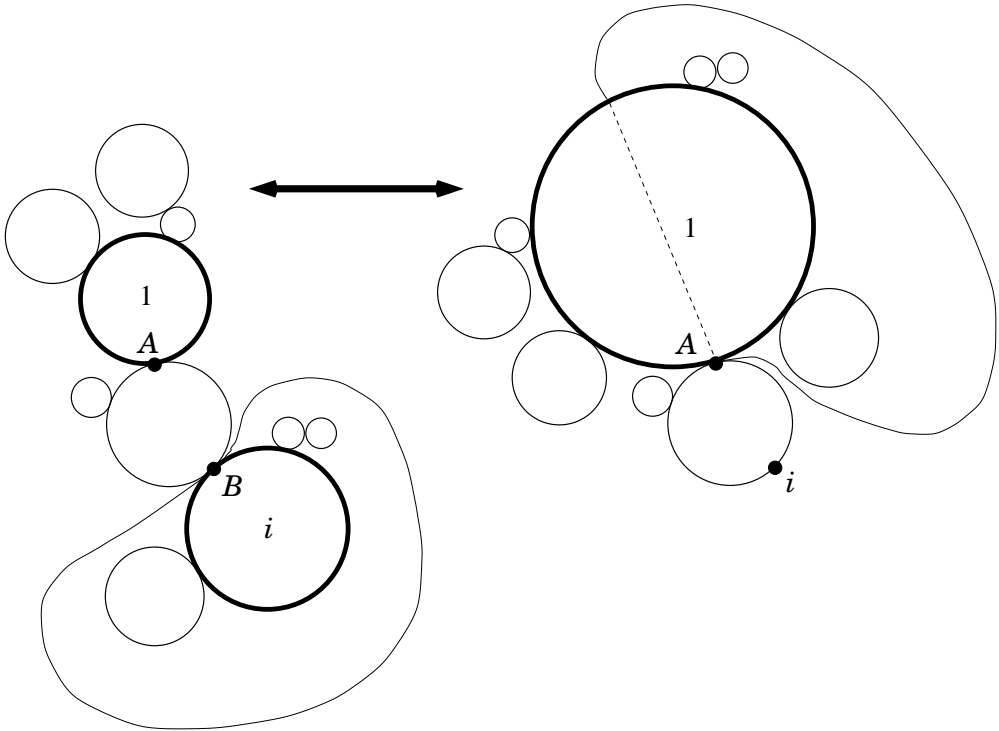}
\caption{Merging two circles.  \label{circlemerge}}
\end{center}
\end{figure}

Suppose we want to merge the circle of color $1$ with the
circle of color~$i$. This merging is shown in 
Figure~\ref{circlemerge}.
 
Since the circle cactus has no cycles, there exists
a unique path of circles leading from the $1$st to the $i$th circle.
Denote by $A$ and $B$ the points of the $1$st and $i$th 
circles respectively, at which this path begins and ends. 
Now we cut our circle cactus in two at the vertex $B$.
The $i$th circle and all circles attached to it belong
to one part, while the other circles belong to the
other part. This cutting is not well defined when more than
two circles meet at the point $B$, but the subspace of 
such circle cacti has zero volume, so it does not matter for us.
We then take the part containing the $i$th
circle and ``open'' it at the point $B$ (so that the $i$th circle
is transformed into an oriented arc). Then we similarly open the 1st circle
at the point $A$, obtaining another oriented arc. The path 
that leads from $A$ to $B$ remains attached to the end (and not
the beginning) of this arc. Finally, we glue the two 
arcs together at their endpoints. The circles
$1$ and $i$ have merged into a single circle.
We mark by the number $i$ the point that was called $B$ before.

It is easy to see that this operation can be reversed
without ambiguity (except on the zero-volume set, where
the point $B$ is a common point of two or more circles).
During the inverse operation the circle of color $1$ 
is cut in two arcs of lengths $l_1$ and $l_i$. Then 
the arcs are closed to form two circles, and, finally,
one of them is moved to the marked point labeled by $i$.

Performing this merging procedure for all $i$ from $2$ to $k$
we end up with a single circle of length $l_1 + \dots + l_k$
with $k-1$ points marked by numbers from $2$ to $k$.
From this marking we can recover
the initial circle cactus, so it can be considered as an
encoding of the circle cactus.

Thus the volume of the space of circle cacti is equal to 
$$
(l_1 + \dots + l_k)^{k-2}.
$$
\cqfd

As a by-product of our considerations we obtain a new
proof of the Cayley formula.

\begin{corollary}
The number of trees with $k$ numbered vertices is equal to
$k^{k-2}$.
\end{corollary}

\paragraph{Proof.} Let us replace every vertex of a tree by a
circle of length 1 and glue circles together iff the corresponding
vertices are joined by an edge. It is easy to see that to
every tree with numbered vertices corresponds a component
of the space of circle cacti whose volume is equal to 1.
Since the total volume of the space of circle cacti
is $k^{k-2}$, there are $k^{k-2}$ components and thus
$k^{k-2}$ trees. \cqfd

\bigskip

Now we will find the volume of the space of circle
cacti in the case when several circles may have the
same color. Suppose we are given $k$ unordered lists
of lengths (a list for each color). Denote by $l_i$
the sum of lengths of circles of color $i$, by
$l$ the sum of lengths of all the circles, and by $m_i$ the
number of circles of color $i$. Further,
denote by $|\Aut_i|$ the number of permutations of 
the circles of color $i$ that preserve their lengths
(for example, if among the circles of color $i$
3 have the same length and all the other lengths
are different, then $|\Aut_i| = 6$).

\begin{theorem} \label{2circ}
In the above notation, the volume of the space of circle
cacti glued from the circles of given lengths is equal to
$$
l^{k-2} \prod_{i=1}^k \frac{(l-l_i)^{m_i-1}}{|\Aut_i|}.
$$
\end{theorem}

\paragraph{Proof.} We are going to apply again the
merging procedure shown in Figure~\ref{circlemerge}, but this time
to circles of the same color. 

The first thing we do,
is to make all the circles of the cactus distinguishable
by labeling (in an arbitrary way) the circles that have the same
color and length. The number of different
labelings of the circles of color $i$ is equal to 
$|\Aut_i|$.

After fulfilling the merging
procedure as many times as possible, we will obtain
a circle cactus with only one circle of every color. Some
points on the circles of this circle cactus will by marked by dots of
different colors in such a way that a dot never lies on the
circle of the same color.

The final circle cactus and the marking allow one to recover the
initial circle cactus without ambiguity.

Finally, the volume of the space of markings of 
any circle cactus with only one
circle of each color (and with given lengths of the circles) 
equals
$$
\prod_{i=1}^k (l-l_i)^{m_i-1},
$$
because the $m_i-1$ marked points of color $i$ lie on the
circles of colors different from $i$, whose total length
is $l-l_i$.

Thus the total volume of the space of circle cacti whose circles
of the same color and length are labeled, is
equal to the product of the volume $l^{k-2}$ 
(given by Theorem~\ref{1circ}) and of the above volume of 
the space of markings.  
In order to obtain the volume of the space of circle cacti
with circles without labeling, we divide this product by the
number
$$
\prod_{i=1}^k |\Aut_i|
$$
of labelings.
The final result is
$$
l^{k-2} \prod_{i=1}^k \frac{(l-l_i)^{m_i-1}}{|\Aut_i|}.
$$
\cqfd

\section{Enumeration of cacti}

The two first theorems of this section were first proved by Goulden and
Jackson (see~\cite{GouJac}) using the Lagrange inversion
theorem; then it was reproved by Lando and the second
author by algebro-geometric methods (see~\cite{LanZvo}).
Here we give elementary proofs. These proofs almost
repeat the proofs of the two parallel theorems
on circular cacti in the previous section (actually,
the copy/paste procedure was used). The difference is
that now we must be careful in counting the cases where
several polygons are glued together at the same vertex,
while in the case of circles such gluings formed a zero-volume
set.

\begin{theorem} \label{1polyg}
Consider the cacti composed of $k$ polygons with, respectively,
$n_1, \dots, n_k$ sides, and whose
colors are all different. Denote 
$$
n = n_1 + \dots + n_k -k+1
$$
be the total number of vertices in such a cactus.
Then the number of such cacti is equal to $n^{k-2}$.
\end{theorem}

\begin{example} \label{ex:3polyg}
Suppose there are 3 polygons with $n_1$, $n_2$, and $n_3$
sides respectively. If the polygon number 1 is in the middle,
the number $\nu$ of edges lying between the 2nd and the 3rd polygons
can vary from $0$ to $n_1-1$ (see Figure~\ref{3polygons}).

\begin{figure}[htb]
\begin{center}
\
\epsfbox{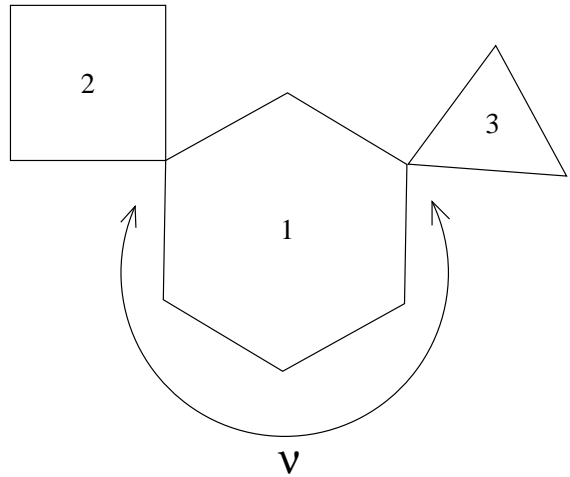}
\caption{A cactus composed of 3 polygons. \label{3polygons}}
\end{center}
\end{figure}

Thus there are $n_1$ possible cacti
with the 1st polygon in the middle. Similarly, there are $n_2$
cacti with the 2nd polygon in the middle and $n_3$ cacti with the
3rd polygon in the middle. But the cactus in which the three 
polygons share a common vertex has been counted three times.
Thus the total number of cacti is $n_1+n_2+n_3 -2 = n$.
\end{example}

\paragraph{Proof of Theorem~\ref{1polyg}.} 
We are going to describe a procedure of
a re-gluing of the cactus that allows one to reduce by $1$ the
number of polygons (and of colors) by merging two polygons into
one. On the new cactus thus obtained we will mark one vertex,
so that the initial cactus can be recovered from the new one.
After applying this merging procedure $k-1$ times, we will
obtain a unique polygon with $n$
vertices, $k-1$ of which are marked. Such polygons will thus
be in a one-to-one correspondence with the cacti, and their
number is equal to $n^{k-2}$.

First of all, we make a list of the numbers $n_i$
that will be needed to reverse the merging operation.

\begin{figure}[htb]
\begin{center}
\
\epsfbox{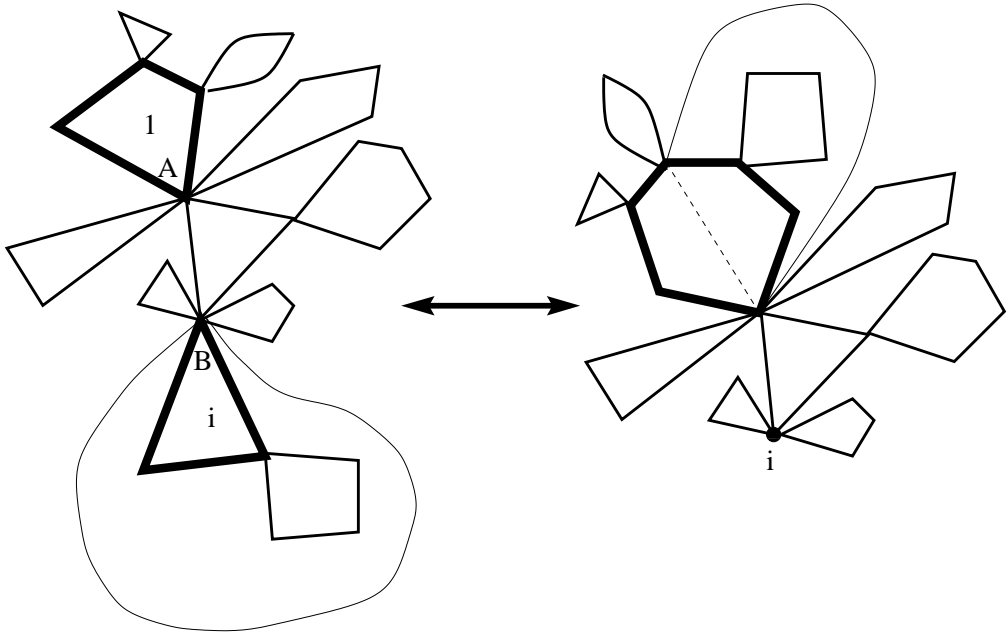}
\caption{Merging two polygons. \label{polygmerge}}
\end{center}
\end{figure}

Let us describe how to merge the polygon of color $1$ with the
polygon of color $i$. This merging is shown in 
Figure~\ref{polygmerge}.
 
Consider the path of polygons leading from the $1$st to the $i$th polygon
(it is unique because the cactus has no cycles).
Denote by $A$ and $B$ the vertices of the $1$st and $i$th 
polygons respectively, belonging to this path. Now we cut
our cactus in two at the vertex $B$.
The $i$th polygon and all polygons attached to it
at vertices different from $B$ belong to one part, while the other
polygons belong to the
other part. We then take the part containing the $i$th
polygon and ``open'' the vertex $B$ (so that the $i$th polygon
is transformed into a broken line). Then we erase the edge of
the $1$st polygon that follows the vertex $A$ in the counterclockwise
direction (thus the $1$st polygon is also transformed into
a broken line). Finally, we glue the two broken lines
together at their free ends. The polygons 
$1$ and $i$ have merged into a single polygon.
We mark by the number $i$ the vertex that was called $B$ before.

It is easy to see that this operation can be reversed
without ambiguity and that it preserves the total number
of vertices. When we reverse the operation, we cut the 
$1$st polygon at one its vertices, so as to obtain
a piece with $n_1-1$ edges and a piece with $n_i$ edges.
The first piece is then closed with an additional edge,
while the second piece is closed by gluing together its
ends. The second polygon is then moved to the marked point
with label~$i$.

Performing this merging procedure for all $i$ from $2$ to $k$
we end up with a single polygon with $n$
vertices, some of which are marked by numbers from $2$
to $k$ (a vertex can be marked by several numbers 
simultaneously). From this marking we can recover
the initial cactus, so it can be considered as an
encoding of the cactus.

Thus the number of such cacti is $n^{k-2}$.
\cqfd

\bigskip

To a $k$-cactus we can assign a list $\X = (X_1, \dots, X_k)$
of $k$ partitions of the number of its vertices.
Their entries are all greater than or equal to 2. The
entries of the partition $X_i$ are simply the numbers of
sides of the polygons colored in the color $i$.

\begin{definition}
The list $\X$ will be called the {\em passport} of the cactus.
\end{definition}

Let $X$ be a partition whose all entries are greater than
or equal to $2$ and suppose that the entry $2$ appears
$a_2$ times, the entry $3$ appears $a_3$ times, etc.
Then we define the number of automorphisms of $X$ to be
$$
|\Aut X| = a_2! a_3! \dots .
$$

Cacti whose polygons are not necessarily of distinct
colors may possess nontrivial symmetries. A {\em symmetry}
of a cactus is a one-to-one map from the sets of its
vertices and edges to themselves preserving adjacency
and colors. Every cactus will be counted with the
weight $1/|\mbox{Sym}|$, where $|\mbox{Sym}|$ is its
number of symmetries.

\begin{theorem}\label{2polyg}
Let $\X = (X_1,\dots, X_k)$ be a passport. Denote
by $p_i$ the number of entries of the partition
$X_i$ and by $n_i$ the sum of these entries. Let
$n$ be the total number of vertices of a cactus with
passport $\X$.

Then the number of such cacti, each cactus
being counted with the weight $1/|\mbox{Sym}|$, is equal to
$$
n^{k-2} \prod_{i=1}^k \frac{(n-n_i)!}{(n-n_i-p_i)! \quad |\Aut X_i|}.
$$
\end{theorem}

\paragraph{Proof.} We are going to apply again the
merging procedure shown in Figure~\ref{polygmerge}, but this time
to polygons of the same color. 

The first thing we do,
is to make all polygons of the cactus distinguishable
by labeling (in an arbitrary way) the polygons that have the same
color and number of sides. The number of different
labelings of the polygons of color $i$ is equal to 
$|\Aut X_i|$.

After fulfilling the merging
procedure as many times as possible, we will obtain
a cactus with only one polygon of every color. Some
vertices of this cactus will by marked by dots of
different colors and these markings will satisfy two
conditions: first, a vertex marked by a dot of color $i$ 
does not belong to the polygon of color $i$; second,
a vertex cannot be marked twice by dots of the same color 
(although it can be marked several times by dots of different
colors).

The final cactus and the marking allow one to recover the
initial cactus without ambiguity.

Finally, any cactus with only one
polygon of each color (and with given numbers of sides
of the polygons) has the same number of possible 
markings equal to  
$$
\prod_{i=1}^k \frac{(n-n_i)!}{(n-n_i-p_i)!} .
$$
On the other hand, the number of cacti whose all polygons
are of different colors is equal to $n^{k-2}$ (Theorem~\ref{1polyg}). 
Multiplying this number by the number of markings, we
obtain the number of cacti with passport $\X$ and with labeled polygons.
In order to obtain the number of cacti
with no labeling of polygons, we divide this product by the
number
$$
\prod_{i=1}^k |\Aut X_i|
$$
of labelings.

We obtain that the number of cacti with passport $\X$
is equal to 
$$
n^{k-2} \prod_{i=1}^k \frac{(n-n_i)!}{(n-n_i-p_i)! \quad |\Aut X_i|}.
$$
\cqfd

\bigskip

Now we are going to consider the case of a rational function
on the Riemann sphere with two poles: one of order $n$ and one
simple. The case of rational functions with $2$ poles
was first studied by V.~Arnold in~\cite{Arnold}.
If we put the poles at $\infty$ and $0$ respectively,
the function will take the form
$$
z^n + a_1 z^{n-1} + \dots + a_0 + \frac{b}{z}
$$
(by a change $z \mapsto az$ the leading coefficient can
be chosen to be equal to~$1$). Counting the number of
rational functions in this case is equivalent to 
finding the multiplicity of the Lyashko-Looijenga map
on the strata of the $D_n$ singularity (see~\cite{Looije}).
We will study only the case where each critical value has only one
critical point (and several non-critical points) as its
preimage. This corresponds to constellations with
polygons of distinct colors. First, let us see what
kind of constellations correspond to the rational functions
of the above form.

\begin{definition}
A {\em $(1,n)$-constellation} is a constellation with
$n+1$ vertices satisfying the following conditions.
Its natural embedding surface is a sphere. The
exterior of the polygons of the constellation on the
sphere consists of two cycles. In one of these cycles,
called the {\em small cycle}, the colors of all the
edges are distinct and their cyclic order is the 
inverse of the cyclic order of the colors. (The small
cycle does not necessarily contain polygons of all
colors.)
\end{definition}

Riemann's existence theorem implies that any $(1,n)$-constellation
determines a rational function of the form
$$
z^n + a_1 z^{n-1} + \dots + a_0 + \frac{b}{z}
$$
uniquely, up to a change of variables
$z \mapsto \varepsilon z$, where $\varepsilon$ is
an $n$th root of unity. The interior of the
small cycle on the sphere contains
the simple pole of the function.

\begin{theorem} \label{3polyg}
Suppose a set of $k$ polygons of
distinct colors is given.
The number of $(1,n)$-constellations glued
from these polygons is equal to
$$
(k-1) n^{k-2}.
$$
\end{theorem}

\paragraph{Proof.} First let us prove that the number
of $(1,n)$-constellations of the above form is equal to
$$
\sum_{
\begin{array}{c}
p \geq 2, \, q \geq 0 \\
p+q \leq k
\end{array}
}
\sum {p+q \choose p} (p-1)^q (m-2p-q) (n+1-p)^{k-p-q-1}.
$$
Here the second sum is taken over all choices of $p+q$
polygons among $k$, and $m$ is the total number of vertices
in these polygons.

Indeed, the number $p$ is
the number of polygons that form the small cycle. The number
$q$ is the number of polygons that do not belong to the
small cycle, but have a vertex lying on it. 
We are going to construct
a $(1,n)$-constellation together with its embedding
into the sphere. First we choose the numbers
$p$ and $q$, then $p+q$ polygons among the $k$.
Then we choose the $p$ polygons that will form the
short cycle, which explains the factor
$p+q \choose p$. There is a unique way to arrange
these $p$ polygons into a cycle respecting the
inverse cyclic order of colors. 

Now we must glue
each of the $q$ polygons to a vertex of the small cycle
respecting the cyclic order of colors at each vertex.
Since each of the $q$ polygons must remain in the exterior
part of the short cycle, there is exactly one vertex of the
short cycle where it cannot be glued. Once we have chosen
the vertices to which the $q$ polygons will be glued,
their cyclic order is automatically determined. Thus 
there are $(p-1)^q$ ways to glue the $q$ polygons.
The obtained graph has $p$ vertices lying on the
small cycle and $m-2p-q$ ``exterior'' vertices.
If, in the total $(1,n)$-constellation, we replace
this subgraph by a unique polygon with $m-2p-q$
vertices, we will obtain a cactus. Conversely,
starting with a cactus containing a $(m-2p-q)$-gon,
we can replace this $(m-2p-q)$-gon with the graph that
we have constructed above and obtain a $(1,n)$-constellation.
Therefore we just need to multiply the number of the
above graphs by the number of cacti given by 
Theorem~\ref{1polyg}:
$$
(n-p+1)^{k-p-q-1}.
$$
This immediately gives the formula
$$
\sum_{
\begin{array}{c}
p \geq 2, \, q \geq 0 \\
p+q \leq k
\end{array}
}
\sum {p+q \choose p} (p-1)^q (m-2p-q) (n-p+1)^{k-p-q-1}.
$$

Now note that the only factor in this formula that depends on
the choice of the $p+q$ polygons (and not only on the numbers
$p$ and $q$) is the factor $m-2p-q$. Therefore we can rewrite the
sum as follows:
$$
\sum_{
\begin{array}{c}
p \geq 2, \, q \geq 0 \\
p+q \leq k
\end{array}
}
{p+q \choose p} (p-1)^q (n-p+1)^{k-p-q-1} \sum (m-2p-q).
$$
Here the second sum is again taken over all the choices
of $p+q$ polygons among the $k$, and $m$ is the total number
of vertices in these polygons. The second sum can now be
evaluated; it is equal to
$$
{k \choose p+q} 
\Bigl[
\frac{(n+k+1)(p+q)}{k} - 2p-q
\Bigr] \,\, .
$$
Indeed, the total number of vertices in all the polygons is equal
to $n+k+1$. Therefore the average number of vertices in $p+q$
polygons is equal to
$$
\frac{(n+k+1)(p+q)}{k} \,\, .
$$

Substituting the value of the second sum in the initial expression,
we obtain that the number of $(1,n)$-constellations is equal to
$$
\frac1{k}
\sum_{
\begin{array}{c}
p \geq 2, \, q \geq 0 \\
p+q \leq k
\end{array}
}
\frac{k!}{p!q!(k-p-q)!} (p-1)^q (n-p+1)^{k-p-q-1}
[np+nq+p+q-kp].
$$
This expression depends only on $n$ and $k$, but not
on the set of polygons that compose the $(1,n)$-constellation.

It can be evaluated by elementary but
heavy computations. First we fix $p$ and carry out
the summation over $q$. Using the formulas
$$
\sum {m \choose r} x^r = (1+x)^m
$$
and
$$
\sum {m \choose r} r x^r = mx(1+x)^m,
$$
one gets
$$
\frac1{nk}
\sum_{p=2}^k 
{k \choose p} \,\,
[np+p-k] \,\,
n^{k-p}.
$$
Now, using once again the same two formulas, the summation
can be carried out over $p$, and it gives
$$
(k-1)n^{k-2}.
$$
\cqfd

\section{Matrix integrals}

In this section we show how to write a matrix integral
that gives a generating function for the volumes of spaces
of circle constellations glued of a given set of colored circles.

\subsection{Matrix models and the Wick formula}

First, let us briefly introduce the physical vocabulary concerning
matrix models and one result that we will need: the Wick formula.
For an accessible introduction to matrix models see~\cite{Zvonkin}.
Multi-matrix models are well explained in~\cite{DiFran}.

Let $E$ be a real vector space (in our case it will be the
space of $k$-tuples of matrices of a special form). It is
called the {\em space of states} of the model. To every point
$x$ of $E$, we assign a nonnegative number $H(x)$, called
the {\em energy} of $x$. (In our case $H$ will be a quadratic form 
on $E$.) The {\em partition function} of the
model is defined by 
$$
Z = \int_E e^{-\frac12 H(x)} dx.
$$
It can be thought of as the number of states of the model, each
state counted with the weight $\exp (-\frac12 \, H)$. For any function $f$
defined on $E$, its {\em expectation
value} is given by
$$
\left< f \right> = \frac1Z \int_E f(x) e^{- \frac12 H(x)} dx.
$$

In the case when $H$ is a positive definite quadratic form on $E$,
the {\em Wick formula} allows one to calculate the expectation
values of polynomial functions on $E$. Denote by $H^{-1}$ the
symmetric bilinear form, dual to $H$, on the dual vector space $E^*$.
(If we choose a basis of $E$ and the dual basis of $E^*$,
then the matrix of $H^{-1}$ is the inverse of the matrix of $H$.) 
Let $\lambda_1, \dots, \lambda_n \in E^*$ be $n$ linear
forms on $E$, and $f = \lambda_1 \dots \lambda_n$ their
product. Then the expectation value of $f$ equals $0$
if $n$ is odd and, if $n$ is even, it is given by the Wick
formula:
$$
\left< f \right> =
\sum_{\mbox{parings of } \{ 1, \dots, n \} } \qquad
\prod_{\mbox{pairs } \{i,j \} } H^{-1}(\lambda_i, \lambda_j),
$$
where a pairing is a way to divide the indices $1,\dots,n$ into
$n/2$ unordered pairs.

Let $E$ be the real vector 
space of $N \times N$ hermitian matrices and
let the energy of a matrix $M$ equal $H(M) = N \tr M^2$. Let
$$
f(M) = \frac{(N \tr M/1)^{k_1}}{k_1!} \,
\frac{(N \tr M^2/2)^{k_2}}{k_2!} \, \dots \,
\frac{(N \tr M^p/p)^{k_p}}{k_p!}.
$$
Then the Wick formula allows one to show that the expectation
value of $f$ equals
$$
\sum \frac1{|\mbox{Sym}|} N^{\chi}.
$$
Here the sum is taken over all not necessarily connected
embedded graphs with $k_1$ vertices of degree~1,
$k_2$ vertices of degree~2, \dots, $k_p$ vertices
of degree~$p$. $|\mbox{Sym}|$ is the number of 
automorphisms of the graph, and $\chi$ is the Euler
characteristic of the surface into which it
is embedded.

The last result and the Wick formula are proved, for 
example, in~\cite{Zvonkin}.

\subsection{Expressing the volume of the space of gluings via
matrix integrals}

Now we will show that the problem of computing
the volume of the space of gluings of circles
can be reduced to a matrix model with $k$ matrices.

We will therefore need a more elaborate version of the
results reviewed in the previous subsection. This version is
fit for colored graphs. The space $E$ will
be the space of $k$-tuples of matrices 
$(A_1, \dots, A_k)$ of a special form 
(a matrix per color). The quadratic form $H$ will encode the
information that only vertices of distinct colors can
be joined by an edge.

According to the Wick formula, the quadratic form that
allows only vertices of different colors to be joined by
edges satisfies
$$
H^{-1} = \sum_{i \not= j} \tr A_i A_j.
$$
The choice of the space $E$ is a consequence of this choice
of $H$. More precisely, consider the complex vector space $E_{\C}$
of $k$-tuples of complex matrices $N \times N$. Then $E$ is
its unique real 
subspace such that (i) $E \otimes_{\R} \C = E_{\C}$ and (ii)
the restriction of $H$ on $E$ is a real positive definite quadratic
form.

Suppose we are given a set of circles colored in $k$ colors
and endowed with lengths. We suppose that the circles of the same color 
and length are labeled, so that all the circles are distinguishable.
Denote by $l_{ij}$ the length of the $j$th circle of color $i$.
We are interested in the following generating function:
$$
F(N,l_{ij}) = \sum_{C} \mbox{Vol} (C) N^{\chi(C)}.
$$
Here the sum is taken over all topological types $C$ 
(see the proof of Theorem~\ref{asymptotic}) of
circle constellations that are not necessarily connected,
but do not have isolated circles. Vol$(C)$ and $\chi(C)$
are, respectively, the volume of the space of circle
constellations of the given topological type, and the
Euler characteristic of the (not necessarily connected)
surface into which the constellation is embedded.
Actually, we take into account only the topological
types such that three circles never share a common point,
because the contribution of the other topological types to
the volume vanishes.

The generating function $F$ is a series in the
variables $N$ and $l_{ij}$. The powers of $N$ that
appear in $F$ are both positive and negative, but
bounded from above. The powers of the $l_{ij}$ are
positive. Suppose a term of $F$ contains
$l_{ij}$ to the power $d-1$. We will soon see that it
takes into account only the volume of those
components of the space of circle constellations, 
on which the $j$th circle of color $i$ touches $d$
other circles. In other words, looking at each term
of $F$ we can say how many vertices each circle
contains. Therefore we can also deduce the number of
faces of the circle constellation.

We will express $F$ using a matrix model with the
following space state $E$ and energy $H$.

Let $E$ be the space of $k$-tuples of matrices
$A_1, \dots A_k$ such that for all $i$, 
$A_i = X + Y_i$, where $X$ is a hermitian matrix and
$Y_i$ are skew-hermitian matrices satisfying
$\sum Y_i = 0$. We identify the space $E$ with its
dual $E^*$ using the symmetric nondegenerate 
(but not positive definite) bilinear form
$$
\left< (A_1, \dots, A_k) | (B_1, \dots, B_k) \right>
 = \tr (A_1B_1) + \dots + \tr (A_k B_k).
$$

More generally, to any symmetric $k \times k$ matrix $S$ we 
can assign a bilinear form on $E$ given by
$$
\sum_{i,j = 1}^k s_{ij} \tr (A_i A_j), 
$$
where $s_{ij}$ are the entries of $S$.

Denote by $I$ the unit $k \times k$ matrix and by $J$
the $k \times k$ matrix whose all entries are equal to $1$.
The energy function $H$ on $E$ is given by the quadratic form
associated to the matrix
$$
\frac{J - (k-1)I}{k-1}.
$$
Thus the dual quadratic form $H^{-1}$ is given by the matrix
$J-I$.

\begin{proposition}
Both quadratic forms $H$ and $H^{-1}$ are positive definite
on~$E$.
\end{proposition}

\paragraph{Proof.} 
The quadratic form associated to the matrix
$J$ is
$$
\sum_{i,j} \tr (X+Y_i) (X+Y_j) =
k^2 \tr X^2 + 2k \tr X (Y_1 + \dots + Y_k) +
\tr (Y_1 + \dots + Y_k)^2
$$
$$
= k^2 \tr X^2,
$$
because $Y_1 + \dots + Y_k = 0$.

The quadratic form associated to the matrix $I$ is
$$
\sum_i \tr (X+Y_i)^2 = k \tr X^2 + 2 \tr X(Y_1 + \dots +Y_k)
+ \sum_i \tr Y_i^2
$$
$$
= k \tr X^2 + \sum_i \tr Y_i^2.
$$
Substituting these formulas in the expressions for $H$ and $H^{-1}$
we obtain
$$
H(A_1, \dots, A_k) = \frac1{k-1}
\left( k \tr X^2 - (k-1) \sum_i \tr Y_i^2 \right);
$$
and
$$
H^{-1}(A_1, \dots, A_k) = k(k-1) \tr X^2 - \sum_i \tr Y_i^2.
$$
Now, for a nonzero hermitian matrix $X$, we have $\tr X^2 > 0$,
while for a nonzero anti-hermitian matrix $Y$, we have
$\tr Y^2 < 0$. It follows that both $H$ and $H^{-1}$ are
positive definite. \cqfd

\bigskip

Now we give an expression of the function $F$ defined
in the beginning of this subsection.
Denote by $d\lambda$ a Lebesgue measure on the vector space $E$.
(It is unique up to a scalar factor and we fix this factor
once and for all.)

\begin{theorem} \label{int}
The generating function $F$ is given by
$$
F(N, l_{ij}) = 
\frac
{
{\displaystyle \int} \, \prod\limits_{i,j} 
\left(N \, \tr
\frac{\displaystyle e^{l_{ij}A_i} - 1}
{l_{ij}}
\right)
\, \, e^{- \frac12  N H} \,\, d\lambda
}
{\int  e^{- \frac12 N H} \,\, d\lambda},
$$
where the integrals are taken over the space $E$.

More precisely, for any positive values of the $l_{ij}$
and for any positive integer $N$, both the integral above
and the series $F$ converge and their values are equal.
\end{theorem}

\begin{figure}[htb]
\begin{center}
\
\epsfbox{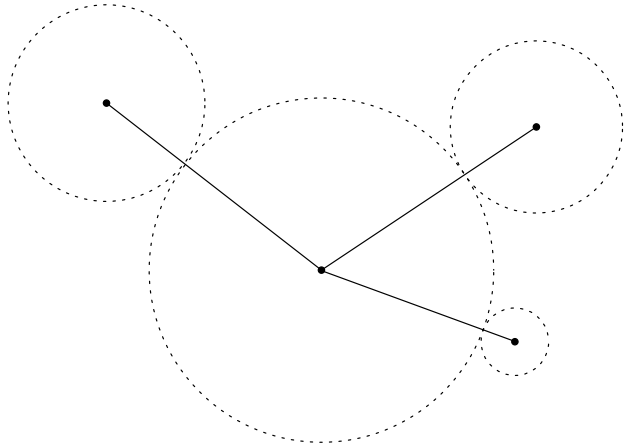}
\caption{Replacing a constellation by its dual. \label{dual}}
\end{center}
\end{figure}

\paragraph{Proof.}
First let us replace every constellation by the
following dual picture (see Figure~\ref{dual}).
Put a vertex in the center of each circle.
Whenever two circles have a common point, draw
an edge trough this common point, connecting the
two corresponding vertices. Thus we
obtain a graph embedded into the surface $\Sigma$. 
This opertion is well-defined whenever no three circles
of the circle constellation share a common point,
i.e., it is well-defined outside a set of measure~0.
The vertices of this graph are colored in $k$
colors, and a number $l_{ij}$ is assigned
to each vertex. Each graph determines
a component of the space of gluings
of circles. Further, for a given gluing
of circles, we can assign a number to each
``angle'' between two edges issued from
the same vertex of the graph. This number
is the length of the part of the circle
that lies between the two corresponding
gluing points.

The sum of numbers surrounding a vertex is equal
to the number $l_{ij}$ assigned to this vertex.
If there are $d$ edges issued from this vertex,
the number $l_{ij}$ must be divided into $d$ 
parts. The volume of the space of such divisions
equals
$$
v_{ij} = \frac{l_{ij}^{d-1}}{(d-1)!}.
$$
If we consider a particular graph $G$, the 
volume of the corresponding component in the 
space of gluings of circles equals
$$
V_G = \prod_{\mbox{vertices of } G} v_{ij}.
$$
The volume of the total space of gluings is
$$
V = \sum_G V_G.
$$

Consider $k$ matrices
$A_1, \dots, A_k$ of size $N\times N$. 
Each matrix corresponds to
vertices of a given color. Consider the
quadratic form 
$$
H^{-1} = \sum_{i \not= j} \tr (A_i A_j)
$$
on the space of the $k$-tuples of matrices.
The form $H^{-1}$ is chosen in such a way, that
when we apply the Wick formula, it would tell us
that we can join
by edges any vertices of different colors, but
not vertices of the same color.

Now we want to enumerate graphs with colored 
vertices and a number $l_{ij}$
attached to each vertex of color $i$. 
If the vertex is of degree
$d$, it will be attributed a weight of
$$
\frac{l_{ij}^{d-1}}{(d-1)!}.
$$
Thus, in the matrix model, such a vertex corresponds
to a term
$$
N \, \frac{l_{ij}^{d-1}}{(d-1)!} \, \frac{\tr A_i^d}d.
$$
Since the degree of each vertex can be
arbitrary, and we want to take into account
all graphs with all possible degrees of vertices,
we must assign to a vertex of color $i$ and
with number $l_{ij}$ the following sum
$$
N \, \sum_{d=1}^{\infty} 
\frac{l_{ij}^{d-1}}{(d-1)!} \, \frac{\tr A_i^d}d
= N \, \tr \frac{e^{l_{ij}A_i} - 1}{l_{ij}}.
$$
Multiplying such factors for all vertices and taking
the expectation value of the product we obtain the
integral from the statement of the theorem.

It remains to prove that 
the operation of taking the infinite sum 
$$
\sum_{d=1}^{\infty} \frac{\tr (l_{ij} A_i)^d}{d!}
$$
commutes with the integration. 
To prove that, note that all the partial finite sums
$$
\sum_{d=1}^D \frac{\tr (l_{ij} A_i)^d}{d!}
$$
are bounded by the same function
$$
N e^{l_{ij}\,|A_i|},
$$
where $|A_i|$ is the sum of absolute values of the coefficients
of the matrix $A_i$. This bound follows from the fact that
the $N$ eigenvalues of $A_i$ are bounded by $|A_i|$. On the
other hand, the product of a finite number of functions
$$
N e^{l_{ij}\,|A_i|}
$$
is integrable with respect to the measure
$$
e^{-N H} d\lambda,
$$
because $H$ is a positive definite quadratic form, while
$|A_i|$ increases only linearly with $A_i$.
It follows that the integral in the statement of the theorem
converges and is equal to the infinite sum of integrals that compute
individual terms of the series $F$. Therefore the series $F$
itself converges to the same value as the integral.
\cqfd

\subsection{Computation of the integral}

We have been able to compute
the integral of Theorem~\ref{int} only
for $N=1$. Therefore the graphs corresponding
to surfaces $\Sigma$ of different Euler
characteristics are enumerated together.

More precisely, suppose as above that we are given
a set of circles colored in $k$ colors and endowed with
lengths $l_{ij}$. Denote by $f(l_{ij})$ the generating
function obtained from $F(N, l_{ij})$ by setting $N=1$.
If we fix the lengths $l_{ij}$ of the circles, the value
of $f$ is the sum of volumes of the spaces of circle
constellations on all compact surfaces. The constellations
(and therefore the surfaces) are not necessarily connected,
but do not contain isolated circles.

In the case $N=1$, the space $E$ is the space of
$k$-tuples of complex numbers $(a_1, \dots, a_k)$
of the form $a_j = x + i y_j$, where $x$ and $y_j$
are real numbers, $\sum y_j = 0$. The quadratic form
$H^{-1}$ becomes
$$
H^{-1} = \sum_{i \not= j} a_i a_j.
$$

\begin{theorem}
$$
f = 
\frac1{\prod\limits_{i,j} l_{ij}} \,\,\,\,
\sum_{{\rm choices \,\,\, of}\,\,\, U_i,\,\, V_i}
(-1)^{\displaystyle |V_1|+ \dots + |V_k|} \,\,\,
e^{\frac12 \, {\displaystyle H^{-1}(L_1, \dots, L_k)}}.
$$
Here the sum is taken over all possible
ways to divide the circles of each color $i$
in two groups $U_i$ and $V_i$. $|V_i|$ is the
number of circles in the group $V_i$. $L_i$
is the sum of lengths of circles in the group $U_i$.
\end{theorem}

\begin{example}
Suppose we have two colors and two 
circles of each color. The lengths
of the circles of color~$1$ are $l_1,l_2$;
those of the circles of color~$2$
are $s_1, s_2$. Then we have
$$
f =  \frac1{l_1l_2s_1s_2}
\,\,
\Bigl[
e^{(l_1+l_2)(s_1+s_2)}
- e^{(l_1+l_2)s_1}
- e^{(l_1+l_2)s_2} 
- e^{l_1(s_1+s_2)}
- e^{l_2(s_1+s_2)} +
$$
$$
\lefteqn
{+ e^{l_1s_1} +  e^{l_1s_2} +  e^{l_2s_1} +  e^{l_2s_2}
-1
\Bigr] =}
$$
$$
= 2 + \frac32 (l_1s_1 + l_1s_2+l_2s_1+l_2s_2) +
\hspace*{10cm}
$$
$$
+ \Bigl[
\frac23 (l_1^2s_1^2 + l_1^2s_2^2+l_2^2s_1^2+l_2^2s_2^2) +
(l_1l_2s_1^2 + l_1l_2s_2^2+l_1^2s_1s_2+l_2^2s_1s_2) + 
\frac32 l_1l_2s_1s_2
\Bigr] + \dots
$$
For instance, the term
$$
\frac23 \, l_1^2 s_1^2
$$
corresponds to the four graphs in Figure~\ref{4graphs}.
In the figure we have drawn the graphs dual to the constellations, 
rather than the constellations themselves. They have two
vertices of color~$1$ with valencies $1$ and~$3$, and
two vertices of color~$2$, also with valencies $1$ and~$3$.

\begin{figure}
\begin{center}
\
\epsfbox{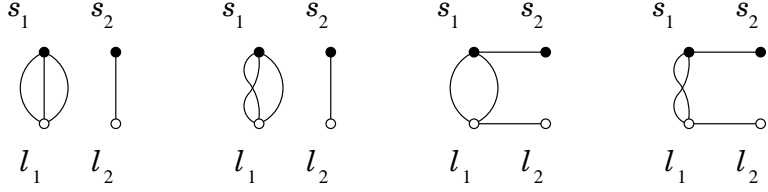}
\caption{The four graphs that contribute to the term
$\frac23 \, l_1^2 s_1^2$. \label{4graphs}}
\end{center}
\end{figure}

Each of these graphs corresponds to a component of
the space of constellations. For the two last graphs,
the volumes of these components equal
$$
\frac{l_1^2}{2!} \, \frac{s_1^2}{2!}.
$$
The first two graphs have an additional symmetry
group of order~$3$, so their components have volume 
$$
\frac 13 \, \frac{l_1^2}{2!} \, \frac{s_1^2}{2!}.
$$
They add up to
$$
\frac{l_1^2 s_1^2}{4} \, (1+1+\frac13 +\frac13)
= \frac23 \, l_1^2 s_1^2.
$$
\end{example}

\paragraph{Proof of the theorem.} 
Let us rewrite the integral to be computed
for $N=1$ :
$$
f = 
\frac
{{\displaystyle \int} \, \prod\limits_{i,j} 
\frac{\displaystyle e^{\displaystyle l_{ij}a_i} - 1}
{l_{ij}}
\, \, 
e^{- \frac12 \,
{\displaystyle H(a_1, \dots, a_k)}} \,\, d\lambda}
{\int  e^{- \frac12 \,
{\displaystyle H(a_1,\dots,a_k)}} \,\, d\lambda}.
$$

\begin{lemma}
\label{fourier}
For any positive definite quadratic form $H$ 
and any linear form $\lambda$ on $E$, we have
$$
\frac
{\int_E e^{- \frac12 H(x) + \lambda(x)}\, dx}
{\int_E e^{- \frac12 H(x)}\, dx}
= e^{\frac12 H^{-1}(\lambda)}.
$$
\end{lemma}

This is a reformulation of the fact that the Fourier
transform of a Gaussian function is the Gaussian function
corresponding the the dual quadratic form. \cqfd

Consider the product
$$
\prod (e^{l_{ij}a_i} -1)
$$
in the numerator of the integral expressing $f$.
Expanding this product we choose either the term
$$
e^{l_{ij}a_i}
$$
or the term $-1$ in each parenthesis. Thus for every $i$,
the circles of colour $i$ are divided in two groups:
we put a circle in the group $U_i$ if we have chosen
$$
e^{l_{ij}a_i}
$$
in the corresponding factor, and in the group $V_i$ if we
have chosen $-1$. Now, the product of any finite number
of terms
$$
e^{l_{ij}a_i}
$$
is the exponent of a linear form on the space $E$.
Therefore for each choice of terms we 
obtain an integral as in Lemma~\ref{fourier}.
Each of them can be calculated, and their sum is given
in the statement of the theorem. \cqfd

\bigskip

\centerline{\bf Acknowledgments}

\medskip

{\small
We are grateful to V.~Arnold, Yu.~Burman, 
M.~Kazaryan, S.~Lando, A.~Zorich, and A.~Zvonkin
for useful and simulating
discussions. We would also like to thank for their
interest the members of
Arnold's seminar in Moscow and of the mathematical
seminars of Rennes University and of Moscow Steklov Institute.
}


\begin{thebibliography}{99}



\bibitem{Arnold} {\bf V.~I.~Arnold.} {\it Topological Classification
of Trigonometric Polynomials and Combinatorics of Graphs with an
Equal Number of Vertices and Edges}, Functional Analysis and its 
Applications,
{\bf vol.~30}, No.~1, 1--17 (1996).


\bibitem{DiFran} {\bf P.~Di~Francesco.} 
{\it Matrix Model Combinatorics: Applications to 
Folding and Coloring.} Random Matrices and their
Applications, MSRI Publications, {\bf vol.~40}, 111-170
(2001).
\verb+http://xxx.lanl.gov+. Archive number: math-ph/9911002.




\bibitem{GouJac} {\bf I.~P.~Goulden, D.~M.~Jackson.}
{\it The combinatorial relationship between trees, cacti
and certain connection coefficients for the symmetric group},
Europ. J. Combinatorics, {\bf vol.~13}, 357--365 (1992).


\bibitem{LanZvo}{\bf S. K. Lando, D. Zvonkine.}
{\em On multiplicities of the Lyashko--Looijenga mapping
on the discriminant strata.}
Functional Analysis and its Applications, 
{\bf vol.~33}, no.~3 (1999)


\bibitem{Looije} {\bf E.~Looijenga.} {\it The complement of the
bifurcation variety of a simple singularity.} Inventiones
Mathematicae, {\bf vol.~23}, 105-116 (1974).


\bibitem{Volklein}   {\bf H. V\"olklein}
{\em Groups as Galois Groups. An Introduction.}
Cambridge University Press,
``Cambridge Studies in Advanced Mathematics'', {\bf vol.~53} (1996).



\bibitem{Zvonkin} {\bf A. K. Zvonkin.}
{\em Matrix integrals and map enumeration: An accessible
introduction.}
Computers and Mathematics with Applications:
Mathematical and Computer Modeling, special issue
``Combinatorics and Physics'' (M.~Bousquet-M\'elou,
D.~Loeb eds.), {\em vol.~26}, no.~8-10, 281--304 (1997).


\end{thebibliography}
\end{document}